\documentstyle[a4,twoside]{article}
\newlength{\longueur}

\newcommand{\motGRAS}[2]
  {\settowidth{\longueur}{#1}#1\hspace{-\longueur}\hspace{#2}#1\hspace{-\longueur}\hspace{#2}#1}

\newcounter{hour}     
\newcounter{hours}    
\newcounter{minute}   
\newcounter{sixty}    

\newcommand{\setclock} {
\setcounter{hour}{\time}
\setcounter{sixty}{60}
\setcounter{minute}{\time}
\divide \value{hour} by \value{sixty}
\setcounter{hours}{\value{hour}}
\multiply \value{hours} by \value{sixty}
\advance \value{minute} by -\value{hours}
}

\newcommand{\clock}{
\the\value{hour}:\ifnum\the\value{minute}<10 0\fi\the\value{minute}
}

\setclock
\newcommand{\ba}{{\bf a}}

\newcommand{\bj}{{\bf j}}
\newcommand{\bk}{{\bf k}}
\newcommand{\bm}{{\bf m}}
\newcommand{\bn}{{\bf n}}
\newcommand{\bt}{{\bf t}}
\newcommand{\bx}{{\bf x}}
\newcommand{\by}{{\bf y}}
\newcommand{\bz}{{\bf 0}}
\newcommand{\bB}{{\bf B}}
\newcommand{\bD}{{\bf D}}
\newcommand{\bX}{{\bf X}}
\newcommand{\SB}[1]{{\bf \mbox{\bf \motGRAS{\scriptsize #1}{0.1mm}}}} 

\def\qed{\hfill\mbox{$\Box$}\medskip}

\def\k#1{\kern#1em}
\def\Ib#1{{I\kern-.25em#1}}
\def\Ibb#1{{I\kern-.23em#1}}

\def\NN{\Ibb N}

\newcommand{\DLZZ}{\settowidth{\longueur}{Z}Z\hspace{-0.9\longueur}Z}

\newtheorem{thm}{Theorem}

\newtheorem{lemma}[thm]{Lemma}
\newtheorem{prop}[thm]{Proposition}

\newcommand{\ord}{\mbox{\rm ord}}

\newcommand{\DX}{\mbox{${\cal DX}$} }
\begin{document}

\begin{titlepage}
\title{Operator Expansion\\ 
in the Derivative and Multiplication by
$x$}
\author{A.\ Di Bucchianico\thanks{Author supported by NATO CRG 930554.}\\
Eindhoven University of Technology\\
Dept.\ of Mathematics \& Computing Science\\
P.O. Box 513\\
5600 MB Eindhoven, The Netherlands\\
{\tt sandro@win.tue.nl}
\and
D.\ E.\ Loeb\thanks{Author partially supported by URA CNRS 1304,
EC grant CHRX-CT93-0400, the PRC Maths-Info, and NATO CRG 930554.}\\
LaBRI\\
Universit\'e de Bordeaux I\\
33405 Talence, France\\
{\tt loeb@labri.u-bordeaux.fr}}
\date{\today, \clock}
\end{titlepage}
\maketitle

\begin{abstract}
We generalize to several variables Kurbanov and Maksimov's result that
all linear polynomial operators can be expressed as a formal sum
$\sum_{k=0}^\infty a_k(X)\, D^k$  
 in
terms of the derivative $D$ (or any degree reducing operator) and
multiplication by $x$. 	In contrast,
we characterize those linear operators that can be expressed as $\sum
_{k=0}^{\infty }f_{k}(D)\, X^{k}$  and give several examples.
Generalizations to several variables and arbitrary degree reducing
operators are considered.
\end{abstract}



\section{Introduction}
In this paper, we will study linear operators on polynomials and give certain 
summation formulas. All
operators are assumed to be linear and their domain and 
range are assumed to be the ring of polynomials $K[x]$ where $K$ is some
field of characteristic zero.

An operator is {\em
shift-invariant} if $QE^a = E^a Q$ for all $a \in K$, where $E^a$ is the shift operator
$E^ap(x) = p(x+a)$. A fundamental result of umbral calculus (see
\cite[Theorem~2, p.~691]{RKO}) is that 
$Q$ is shift-invariant if and only if it can be expressed (as an
infinite series) in the derivative $D$. Its coefficients are given by
the D-expansion formula
\begin{equation}\label{expansion}
Q=\sum_{k=0}^\infty a_k D^k 
\end{equation}
 where $a_k = [Qx^k/k!]_{x=0}.$

It has been asked \cite[Problem 12, p.\ 752]{RKO} what operators can be expressed
(as infinite  
series) in $D$ and in $X$ where $X$ represents
multiplication by $x$. Pincherle and Amaldi \cite{PiA}
 proved
that all linear 
operators can be so expressed. Kurbanov and Maksimov give an explicit
construction of such an expansion in
the Proceedings of the Uzbekistan Academy of Science \cite{Kurb} 
(even with the added
restriction that the derivatives lie to the right of multiplication by $x$)
\begin{equation}\label{xD}
Q = \sum_{n= 0}^\infty a_{n}(X)\, D^n
\end{equation}
where the polynomials $a_n$ are defined by the following generating function
\begin{equation}\label{star}
\frac{Q\exp(xt)}{\exp(xt)} = \sum_{n\geq 0} a_n(x)\, t^n .
\end{equation}
(By abuse of notation, $Q$ above acts on the coefficients of
$\exp(xt)\in (K[x])[[t]]$.) 
Such an expansion will be called an
XD-expansion since it is an expansion in $D$ and $X$
with $X$ lying to the left of $D$. Note that XD-expansions are useful in that the
derivative is easy 
to calculate both numerically and symbolically. Moreover, XD-expansions allow us
to manipulate arbitrary linear operators with 
similar ease. For example, the 
XD-expansion of the integration operator gives an asymptotic formula
for integration (\ref{J}) which can also be found in \cite[Sections
3.5 and 3.6]{NB}. Note that infinite sums such as (\ref{xD}) are
always well defined since when applied to any given polynomial $p(x)$
only a finite number of terms make a nonzero contribution.

Equation (\ref{xD}) can be viewed as the XD-generalization of
equation (\ref{expansion}). In fact, many of the results of umbral
calculus are direct corollaries of (\ref{xD}). In section \ref{xds}, a generalization 
of equation
(\ref{xD}) will be proven. As do Kurbanov and Maksimov,  the
derivative will be replaced by a general degree reducing operator $B$. We
will further generalize equation (\ref{xD}) and consider 
multivariate polynomial operators.
 
The duality between the operators $D$ and $X$ is a recurrent theme in 
the umbral calculus
(see for example \cite[p.~695]{RKO} and \cite{Romanbook}). Thus, it is
natural to consider the dual 
form of (\ref{xD}), namely, expansions of the form
\begin{equation}\label{Dx}
Q = \sum_{n=0}^\infty a_{n}(D)\, X^{n}
\end{equation}
called DX-expansions.
In \cite[theorem 2.1]{GJ}, 
we find DX-expansions of certain umbral
operators. Do all 
operators have a DX-expansion?
Surprisingly, the answer is ``no'', and those that do will be called
DX-operators.  

In section \ref{xds}, we will show that the set of DX-operators form a
subalgebra of the algebra of 
linear operators. We will characterize DX-operators in several ways.

There are several surprising differences between DX and XD-expansions
that will be explored in this paper:
\begin{itemize}
\item Although all operators have
XD-expansions, not all have DX-expansions. 
\item Although all XD-expansions converge 
when applied to polynomials, 
not all DX-expansions do so. 
\item Although Maksimov and Kurbanov's formula (\ref{star})
generalizes well to
multivariate polynomial operators, and arbitrary degree reducing
operators $B$, the corresponding  results concerning DX-expansions do
not hold in the more general setting.
\end{itemize}
DX and XD-expansions for many important operators
are given in section \ref{s3}. Examples
include several new operator expansions 
as applications to the umbral calculus; in
particular, we will show that all umbral operators, and umbral shifts
are DX-operators. 

Some open problems are suggested in section \ref{op}.

\section{$XD$-Expansion}\label{xd}\label{xds}

The objective of this section is to prove the XD-expansion formula
(\ref{xD}). However, we will 
first restate the formula in greater generality. Since equation
(\ref{xD}) is a generalization of the D-expansion formula
(\ref{expansion}), we will as an introduction
indicate the corresponding generalization of equation (\ref{expansion}).

Let $\bx=(x_{1},x_{2},\ldots)$ be a finite or infinite set of
variables.\footnote{Note that the reader
interested only in the univariate case may skip all the definitions
involving bold-faced symbols with the assurance that in the univariate
case they reduce to their univariate counterparts. In particular,
$e_{i}$ below 
should then be read as ``1.''} A {\em monomial} over $\bx$ is a product
$\bx^{\SB{n}} = x_{1}^{n_{1}}x_{2}^{n_{2}}\cdots $. Here and in the sequel
$\bn$, $\bm$, $\bj$, $\dots$ denote  sequences of nonnegative integers with 
finite support. A
{\em polynomial} is a finite $K$-linear combination of monomials.
We denote the ring of polynomials over $\bx$ by $K[\bx]$.  

Let $D_{i}$ be the derivative with respect to $x_{i}$. Thus,
$\bD^{\SB{k}}=D_{1}^{k_{1}}D_{2}^{k_{2}}\cdots $. Define 
$\bn!=n_{1}!\,n_{2}!\,\cdots $, $(\bn)_{\SB{k}}=\bn !/ (\bn -\bk)!$, and
${\SB{n}\choose \SB{k}} 
=(\bn )_{\SB{k}}/ (\bn-\bk)!$ with the motivation
that  
$\bD^{\SB{k} }\bx^{\SB{n}} = (\bn)_{\SB{k}} \bx^{\SB{n} - \SB{k}} $ and 
$(\bx + \by)^{\SB{n} } = \sum _{\SB{k}} {\SB{n} \choose \SB{k}}
\bx^{\SB{k}} \by^{\SB{n} - \SB{k}}. $

An operator is shift-invariant if $QE^{\SB{a}} = E^{\SB{a}} Q$  where
$E^{\SB{a}} = \exp(\ba \cdot 
\bD)$ is the
shift operator $E^{\SB{a}} p(\bx) = p(\bx + \ba)$.  We can now generalize
the D-expansion formula (\ref{expansion}).
\begin{prop}[$\bD$-Expansion Formula] A linear operator
$Q \colon 
K[\bx]\rightarrow K[\bx]$ is shift-invariant if and only if it can
be expressed (as an infinite series) in the derivative \bD
$$ Q = \sum_{\SB{k}} a_{\SB{k}} \bD^{\SB{k}}$$	
where $a_{\SB{k}} = [Q\bx^{\SB{k}} / \bk!]_{\SB{x} = \SB{0}}.$
\end{prop}
Another generalization of the D-expansion formula (\ref{expansion}) is
to replace the derivative $D$ with an arbitrary degree reducing
operator $B$. Suppose $B \colon 
K[x]\rightarrow K[x]$ is such that for all
nonconstant polynomials $p$, $(\deg{Bp})+1=\deg{p}$, and 
for all constants $a$, $Ba=0$. Then $B$ is called {\em degree
reducing} and there exists a unique sequence of
polynomials $b_n(x)$ (called the {\em divided power sequence for $B$} 
such that $Bb_n(x)= b_{n-1}(x)$ and
$b_n(0)=\delta_{n0}$. For example, if $B=D$,  then $b_{n}(x)=x^{n}/n!$.
The sequence $n!\,b_{n}(x)$ is called the {\em basic family} of $B$ by Markowsky
\cite{Mark}.

We then have the ``B-expansion formula.''
\begin{prop}[$B$-Expansion Formula] A linear operator
$Q \colon 
K[x]\rightarrow K[x]$ commutes with $B$ if and only if it can
be expressed (as an infinite series) in $B$.
$$ Q = \sum_{k=0}^\infty a_k B^k$$
where $a_k = [Q b_k (x)]_{x=0}$.
\end{prop}
It is of course possible to carry out both generalizations
simultaneously. Let $\{b_{\SB{n}} (\bx) \colon 
\bn \mbox{ finite support}
\}$ be a basis of $K[\bx ]$ such that $b_{\SB{n}} (\bz ) =
\delta_{\SB{n},\SB{0}}$. 
Define the operator $B_{i}$ by  
$B_{i}b_{\SB{n}}(\bx ) = b_{\SB{n} - e_{i}}(\bx )$ where $e_{i}$ is the
$i$th unit vector $(e_{i})_{j}=\delta_{ij}$. Denote by \bB\ the sequence of
operators $(B_{1},B_{2},\ldots)$. For example, if $b_{\SB{n}}(\bx )=
\bx ^{\SB{n}}/\bn !$, then $\bB = \bD. $
\begin{prop}[$\bB$-Expansion Formula] A linear operator $Q \colon
K[\bx ]
\rightarrow K[\bx ]$ commutes with all $B_{i}$ if and only if it can
be expressed (as a formal power series) in \bB.
$$ Q=\sum _{\SB{k}} a_{\SB{k}} \bB^{\SB{k}} $$ 
where $a_{\SB{k}}=[Qb_{\SB{k}}(\bx )]_{\SB{x} =\SB{0}}.$\qed 
\end{prop}
A formal power series is an infinite linear combination of monomials. 
The ring of formal power series with coefficients in $K$ is denoted
$K[[\bx ]].$ The generating function $b(\bx ,\bt)$ for the sequence of
polynomials $b_{\SB{n}}(\bx )$ is 
$$ b(\bx ,\bt) = \sum _{\SB{n}} b_{\SB{n}}(\bx )\, t^{\SB{n}} \in (K[\bx
])[[\bt]].$$
For example, if $\bB = \bD $, then $b(\bx , \bt) = \exp(\bx \cdot \bt).$

By abuse of notation, we allow operators on $K[\bx ]$ to act on the
coefficients of $b(\bx ,\bt).$ For example,
\begin{eqnarray*}
B_{i}\,b(\bx ,\bt) &=& \sum _{\SB{n}} (B_{i}b_{\SB{n}}(\bx ))\, t^{\SB{n}} \\
&=& \sum _{\SB{n}} b_{\SB{n} - e_{i}}(\bx )\, t^{\SB{n}} \\
&=& \sum _{\SB{m}} b_{\SB{m}}(\bx )\, t^{\SB{m} + e_{i} } \\
&=& t_{i}\, b(\bx , \bt).
\end{eqnarray*}
Let $X_{i}$ be the operator of multiplication by $x_{i}$, and let
$\bX=(X_{1},X_{2},\ldots)$. 
All linear operators can be expressed in terms of $\bX$ and $\bB$. 
\begin{thm}[$\bX\bB$-Expansion Formula]\label{XBExp}
Let $\bB=(B_{1},B_{2},\ldots)$, and $b(\bx ,\bt)$ be as above. 
Let $Q \colon 
K[\bx ]\rightarrow K[\bx ]$ be a linear operator. 
Then 
\begin{equation}\label{XB}
Q=\sum _{\SB{k}} a_{\SB{k}}(\bX)\, \bB^{\SB{k}}
\end{equation}
where the polynomials $a_{\SB{k}}$ are 
 given by the generating function
$$ \sum_{\SB{k}} 
a_{\SB{k}} (\bx )\, t^{\SB{k}} = \frac{Qb(\bx ,\bt)}{
b(\bx ,\bt)}.$$ 
\end{thm}

{\em Proof:} Apply both sides of equation (\ref{XB})
to the basis $b_{\SB{k} }(\bx )$, or equivalently to its generating
function $b(\bx ,\bt).$ The 
left-hand side gives $\bB b(x,t)$ while the right-hand side gives
\begin{eqnarray*}
\sum _{\SB{k}} a_{\SB{k}}(\bX)\, \bB^{\SB{k}} b(\bx \cdot \bt) &=& 
\sum _{\SB{k}} a_{\SB{k}}(\bx)\, \bt^{\SB{k}}\, b(\bx \cdot \bt) \\
&=& b(\bx \cdot \bt)^{-1}\,(Qb(\bx \cdot \bt))\, b(\bx \cdot \bt)\\
&=& Qb(\bx \cdot \bt).
\end{eqnarray*}
\qed 
\begin{prop}[$\bX\bB$-Uniqueness] The $\bX\bB$-expansion 
given in theorem \ref{XBExp} is unique.
\end{prop}
{\em Proof:} It  suffices to show that $\sum _{\SB{n}}
a_{\SB{n}}(\bX)\, \bB^{\SB{n}}$ is the zero operator only if
$a_{\SB{n}}(\bX)$ is 
zero for all $\bn$. Suppose not, and let $\bm $ be a minimal such $\bn$.
However, $\left(\sum _{\SB{n}} a_{\SB{n}}(\bX)\, \bB^{\SB{n}}
\right)\bx^{\SB{m}}=\bm!\, a_{\SB{m}}(\bx) \neq 0.$\qed 

We reserve most of our examples for section \ref{s3}. However, we can
not resist giving a simple but important example here.

Let $J$ denote the definite integral $Jp(x) = \int_{0}^{x} p(u)du$.
To express $J$ in terms of $D$ and $X$, we apply $J$ to $\exp(xt)$ 
which yields $(\exp(xt)-1)/t.$ We then divide by $\exp(xt)$ and replace $x$
and $t$ with $X$ and $D$ respectively, being sure to keep $X$ on the
left and $D$ on the right since they do not commute. Thus,
\begin{equation}\label{J}
J= \sum_{n=0}^{\infty }(-1)^{n}X^{n+1}D^{n}/(n+1)!. 
\end{equation}
This single equation is essentially equivalent to the content of
\cite[sections 
3.5 and 3.6]{NB}. 

However, note that the methods used here are purely formal. Thus, with
equal ease, any linear operator in any number of dimensions can be
expanded in terms of the derivative or any other ``delta set'' of
degree reducing operators.
For example, $J$ can be expanded in terms of $X $ and $B=\Delta$ instead
of $X$ and $D$ where $\Delta$ is the forward difference operator
$\Delta p(x) = p(x+1)-p(x).$  The divided power sequence for $B=\Delta$
is $b_{n}(x)={x\choose n}$ whose generating function is given by 
$b(x,t)=(1+t)^{x}$. Thus, $J=\sum _{n=0}^{\infty } a_{n}(X) D^{n}$
where $a_{n}$ is given by the generating function $\sum _{n=0}^{\infty
} a_{n}(x)t^{n}= {1-(1+t)^{-x}}/{\ln(1+t)} $.  
In other words,
$$J=X-X^{2}\Delta/2 +
(X^{2}/4+X^{3}/6)\Delta^{2}-(X^{2}/6+X^{3}/6-X^{4}/24)\Delta^{3}+\cdots.
$$  
Recall that differentiation is simpler than 
integration, and that finite differences are even simpler still. Thus,
formulas such as those above give simple algorithms by which one can
approximate complicated linear operations. For more examples, see \cite{Kurb}. 

\begin{prop}[$\bX\bB$-Convergence] \label{xdconv}
Any infinite sum of the form $\sum _{\SB{n}} a_{\SB{n}} (\bX )
B^{\SB{n}}$ converges formally. In other words, only finitely many terms
are nonzero when applied to any given polynomial.
\end{prop}

{\em Proof:}  Every polynomial
can be expressed as a linear combination
$$ p(\bx ) = \sum _{\SB{n} 
 \in S} c_{\SB{n}} b_{\SB{n}}(\bx )$$ 
for some finite set $S$. Let $m_{i}=\max _{\SB{n}\in S}(n_{i})$. Then
$c_{\SB{n}}(\bX) \bB^{\SB{n}} p(\bx) $ is zero except when $\bn \leq
\bm$ componentwise,
and there are only finitely many such $\bn . $\qed 

Of course, the results above have only been proven for polynomials. One would
hope to extend the expansion formulas by continuity to functions that are
limits of polynomials in some sense. For example, $J \cos x = \sin
x$, 
 and   
\begin{eqnarray*}
{\sum_{n=0}^{\infty}\frac{(-1)^{n}X^{n+1}D^{n}}{n+1} \cos x}
& = &
\sum_{k=0}^{\infty} \frac{x^{2k+1}}{(2k+1)!} (-1)^k \cos x + 
\sum_{k=0}^{\infty} \frac{x^{2k+2}}{(2k+2)!} (-1)^{k} \sin x \\
&=& \sin x \cos x + (1 - \cos x)\sin x\\
& =& \sin x.
\end{eqnarray*}
This is an example of a direct verification. In general, one has to
justify 
the interchange of a limit and a summation. 
The following proposition gives a sufficient condition to allow this
interchange.
\begin{prop}\label{analytic} Let $Q$ be a shift-operator with 
representation $Q = \sum_n a_{n}(X)\, D^n$ and let $f$ be the pointwise
limit of a sequence of polynomials $(p_k)_{k \geq 0}$ for all $x$ in some set $V$. 
 If there exists $\varphi_{n}(x)$ such that $|a_{n}(x) D^{n} p_{k}(x)|
\leq \varphi_{n}(x)$ and
$\sum_{n=0}^{\infty }
 \varphi_{n}(x) < \infty$ 
for all $n,k\in \NN $ and $x \in V$, then 
$$\lim_{k \rightarrow
\infty} \sum_{n=0}^{\infty} a_{n}(X)\, D^n p_k(x) =
\sum_{n=0}^{\infty} a_{n}(X)\, D^n f(x).$$
\end{prop}
{\em Proof:} This follows from Lebesgue's dominated convergence
theorem, since series are integrals with respect to the measure $\sum_n
\delta_n$, where $\delta_n$ is the Dirac measure at $n$.\qed  
\section{$DX$-Expansions}\label{dxs}
\subsection{Commutation Rules}
When an operator has a DX-expansion, the coefficient of $D^{k}X^{n}$ in
(\ref{xD}) is not necessarily the same as that of $X^{n}D^{k}$ in
(\ref{Dx}) since $D$ and $X$ do not commute (cf. proposition
\ref{coeffs}).  In fact, 
their commutator $DX-XD$ is the identity. Applying theorem~\ref{XB} , we obtain
\begin{equation}\label{two}
D^j X^i = \sum_{k=0}^j (i)_k (j)_{k} X^{i-k} D^{j-k} / k!
\end{equation} 
where the lower factorial $(x)_n = x(x-1)(x-2)\dots (x-n+1)$. For another
proof, see \cite[Section~2.2]{FS}.

Let $a_{nm}=D^{n}X^{m}/n!m!$ and 
$b_{nm}=X^{m}D^{n}/n!m!$. Then we have $a_{nm}=\sum _{t}
{b_{n-t,m-t}}/t!$,  
or fixing the difference between $n$ and $m$ and indexing only by $n$,
we have $a_{n}=\sum _{t} b_{n-t}/t!$. This sum is easily inverted, for
the vector ${\bf a}=M{\bf b}$ 
where $M$ is the upper triangular infinite Toeplitz matrix with $i$th
diagonal $1/i!$. However, the inverse of this matrix has diagonals
$(-1)^{i}/i!$ as can be seen from the identity $\sum
_{i+j=k}(-1)^{i}/i!j!= (1-1)^{k}/k! = \delta_{k0}$.
Thus, we obtain
\begin{equation}\label{one}
X^i D^j = \sum_{k=0}^j (-1)^k (i)_k (j)_{k} D^{j-k} X^{i-k} / k! .
\end{equation}

More generally, we have the following result. 
\begin{prop}[${[D,X]}$-Commutation]
\begin{eqnarray}
\label{7.5}
f(D)p(X) &=& \sum _{k=0}^{\infty} p^{(k)}(X)\, f^{(k)}(D) / k!\\
\quad
p(X)f(D) &=& \sum _{k=0}^{\infty} (-1)^{k} f^{(k)}(D)\, p^{(k)}(X) / k! .
\end{eqnarray}
\end{prop}
{\em Proof:} Consider linear combinations of (\ref{two}) and
(\ref{one}).\qed 

In so far as possible, we will derive a ``B'' analog and a
multivariate analog of each result. These may be omitted by the reader
essentially interested in the DX-expansions.

In this case, an arbitrary operator $B$ may have virtually any
commutator $BX-XB$. 
Thus, no corresponding ``[B,X]-commutation'' result
can be given. 
However, in the multivariate case, we obtain the following result.
\begin{prop}[${[\bD ,\bX]}$-Commutation]
\begin{eqnarray*}
f(\bD)p(\bX)\, &=& \sum _{\SB{k}} p^{(\SB{k})}(\bX)\, f^{(\SB{k})}(\bD) /
\bk! \\
p(\bX)\, f(\bD) &=& \sum _{\SB{k}} (-1)^{|\SB{k}|}\, f^{(\SB{k})}(\bD)\,
p^{(\SB{k})}(\bX) / \bk! 
\end{eqnarray*}
where $|\bk  |=k_{1}+k_{2}+\cdots $  and $f^{(\SB{k})}(\bD )$ is given
by the multivariate Pincherle derivative $f^{(\SB{k})}(\bx ) =
D^{\SB{k}}f(\bx )$.  
\end{prop}
{\em Proof:} Apply $[D,X]$-commutation to each variable. Recall that
$D_{i}$ and $x_{j}$ commute for $i\neq j.$\qed 

\subsection{Convergence of $DX$-Expansions}

In section \ref{xd}, we have seen that all operators have
XD-expansions, and even XB, \bX\bD, \bB\bD-expansions.
One might suspect that the commutation rules of the previous subsection 
would lead to corresponding DX ({\em et al})-expansions.
Nevertheless, we will see below that certain operators (for example $J$)
lack a corresponding DX-expansion. (See remark after theorem
\ref{p3}.) The finite commutation 
rules of the previous subsection do not generalize to infinite commutation 
rules, for 
if we substitute (\ref{one}) into (\ref{xD}), the
resulting sum is not necessarily well defined. For example, each
coefficient of $X^{n}D^{n}$ in an XD-expansion
makes its own contribution to the coefficient of $X^{0}D^{0}$ in the
corresponding DX-expansion. Since there are no 
conditions on the coefficients of an XD-expansion, the sum obtained by
direct substitution need not 
converge. That is to say,
when an operator is applied to a polynomial of
degree $n$, the sum in (\ref{xD}) is well defined since only the first
$n+1$ terms make nonzero contributions. On the other hand, all terms
of (\ref{Dx}) make potential contributions. 

What DX-expansions converge? 

To answer this question, we define the {\em order} $\ord (f)$ of
a nonzero formal power series $f(t)=\sum _{k=0}^{\infty } c_{k}t^{k}$ 
to be the smallest $j$ such that $c_{j}\neq 0$.
If $\deg(p)=n$ and $\ord
(f)=k$, then $f(D)p(x)$ 
is of degree $n-k$ if $n\geq k$, and $f(D)p(x)=0$ otherwise. The
order of the zero series is taken to be $+\infty$. 
\begin{prop}[$DX$-Convergence] \label{dxconv}
{\bf (1)} The formal sum of operators $\sum_{k=0}^\infty  f_{k}(D)\, X^k$
converges in  
the discrete topology if and only if $\displaystyle \lim_{k\rightarrow
\infty } 
\left[ \ord 
(f_{k}) - k \right] = +\infty $.

{\bf (2)}The formal sum of operators $\sum_{k=0}^\infty  D^k a_k(X)$
converges in  
the discrete topology if and only if $\displaystyle \lim_{k\rightarrow
\infty } \left[ 
k-\deg(a_{k})  \right] =  +\infty $.
\end{prop}

{\em Proof:} We will prove only {\bf (1)} since {\bf (2)} follows by
similar reasoning.

{\bf (Only if)}
By hypothesis,  the polynomial
sequence  $\left( \sum_{k=0}^K f_{k}(D)\, X^k p(x)\right)_{K\geq 0}$ is
eventually  constant 
for all polynomials $p(x)\in K[x]$. 
Thus, $\sum_{k=K}^{\infty } f_{k}(D)\, x^k p(x)$ is eventually zero.
Hence, $\ord(f_k)-k$ is eventually greater than $\deg(p)$. Since
$p(x)$ may have any degree, $\lim_{k\rightarrow \infty} \ord(f_k)-k
= +\infty$.

{\bf (If)} Let $p(x)$ be a polynomial of degree $n$. Since  
$\lim_{k\rightarrow \infty} \ord(f_k)-k
= +\infty$, there exists a $k_{0}$ beyond which $\ord (f_{k}) - k > n$.
Thus, the sequence $\left( \sum_{k=0}^K f_{k}(D)
X^k\right)_{K \geq 0}$ is constant for 
$K>k_{0}.$\qed 

The ``$B$-analog'' of proposition \ref{dxconv} is very easy to state.
In fact, its proof identical to that of proposition \ref{dxconv} {\em
mutatis mutandis}.
\begin{prop}[$BX$-Convergence]
{\bf (1)} The formal sum of operators $\sum_{k=0}^\infty  f_{k}(B)\, X^k$
converges in  
the discrete topology if and only if $\displaystyle \lim_{k\rightarrow
\infty } \left[ \ord  (f_{k}) - k \right] =  +\infty $.

{\bf (2)} The formal sum of operators $\sum_{k=0}^\infty  B^k a_k(X)$
converges in  the discrete topology if and only if $\displaystyle
\lim_{k\rightarrow \infty } \left[ 
k-\deg(a_{k})  \right] = +\infty $.
\end{prop}
The obvious corresponding result is false in the multivariate case. 
Consider for example the sum $\sum _{k=0}^{\infty }
D_{1}^{k}X_{2}^{k}$ which does converge formally.
The correct necessary and sufficient conditions are more complex than
in the univariate case.
\begin{prop}[$\bD \bX$-Convergence] \label{DXConv}
The operator sum $Q=\sum _{\SB{k}, \SB{n}} c_{\SB{k}, \SB{n}}
\bD^{\SB{k}}  \bX^{\SB{n}}$
converges formally if and only if  
for all \bj, there exist only finitely many triples \bm, \bk,  \bn
such that $\bm + \bn - \bk = \bj $   and
$c_{\SB{k},\SB{n}} \neq 0$. 
\end{prop}
{\em Proof:} The sum converges formally if and only if when applied to any
polynomial there are only finitely many contributions. In other words,
$Q\bx^{\SB{m}}$ has only contributions to finitely many terms
$\bx^{\SB{j}}$, and  only 
finitely many contributions to each such term. These two constraints
correspond exactly to the condition above.\qed 

We will not state the obvious  \bB \bX-analog of proposition
\ref{DXConv}.  

\subsection{Characterization of $DX$-Operators}

What operators can be represented by DX-expansions? 

A simple criterion is given in theorem \ref{p3} based on the
matricial representation of an operator.
Any operator can be represented as an infinite matrix
$C=(c_{nk})_{n,k}\geq 0$ with respect to, for
example, the basis $\{x^{n}\colon n \geq 0\}$,   
$$Qx^n=\sum_{k=0}^\infty c_{nk} x^k .$$
Let $q_{t}$ be the $t$-th diagonal of the matrix $C$:
$q_{t}(n)=c_{n,n+t}$ where $n$ is a nonnegative integer, and $t$  
is an arbitrary integer under the convention that $c_{nk}=0$ for
$k<0$. 

\begin{thm}[$DX$-Characterization]\label{p3}
Using the above notation, 
$Q$ has a DX-expansion if and only if $q_{t}\in K[n]$ for all integers
$t$. 
\end{thm}

While one can see that an operator has a DX-expansion simply by
exhibiting such an expansion, the only convenient way to see that an
operator does not have such an expansion is to apply the criterion of
theorem \ref{p3}. For example, consider the definite integration
operator $J$ defined above. For $Q=J$, $Qx^{n}=x^{n+1}/(n+1)$. Thus,
$q_{1}(n)=1/(n+1)\not\in K[n]$, so $J$ does not have a DX-expansion. 
(See propositions \ref{p5} and \ref{p7} for other examples.)

Theorem \ref{p3} is all the more surprising since it has no
obvious ``BX-analogue.'' For example, let $Bx^{n}=\sin(n)x^{n-1}$.
Now, $Q=B$ 
obviously has a BX-expansion, namely $B$ itself, yet
$q_{-1}(n)=\sin(n) \not\in K[n].$  
On the other hand, theorem \ref{p3} has the multivariate analog
below which we will prove in the place of theorem \ref{p3}.
\begin{thm}[$\bD \bX $-Characterization]\label{bolddxchar}
Let $Q \colon 
K[\bx] \rightarrow K[\bx]$ be a linear operator. Let 
$q_{\SB{t}}(\bn )=c_{\SB{n} ,\SB{n} +\SB{t}}$ where
$Q\bx^{\SB{n} }=\sum _{\SB{k}} c_{\SB{n},\SB{k}} \bx^{\SB{k}}$.
Then $Q$ has a \bD \bX-expansion
$Q=\sum _{\SB{n}} a_{\SB{n}}(\bD)\, (\bx)^{\SB{n}}$
if and only if
$q_{\SB{t}}(\bn )\in K[\bn ]$ for all integer vectors $\bt$ (with
finite support). 
\end{thm}

{\em Proof:} {\bf (Only if)} The only terms which contribute to
$q_{\SB{t}}$ are $\bD^{\SB{k}}\bX^{\SB{t}+\SB{k}}$. By (1) of
proposition \ref{DXConv}, there are finitely many such terms.
Each makes a contribution of $(\bn + \bt + \bk)_{\SB{k}}\in K[\bn ]$. 
Thus, $q_{\SB{t}}(\bn )\in K[\bn ].$

{\bf (If)} Define the operator $Q_{\SB{t}}$ by the relation
$Q_{\SB{t}}\bx^{\SB{n}}= 
q_{\SB{t}}(\bn)\bx^{\SB{n}+\SB{t}}$. Then $Q=\sum _{\SB{t}}Q_{\SB{t}}$
is a convergent 
expansion of $Q$ since $q_{\SB{t}}$ has finite support,
so it will suffice to show that 
each $Q_{\SB{t}}$ is a $\bD\bX$-operator.

The polynomial $p_{\SB{k}}(\bn) = (\bn+\bt+\bk )_{\SB{k}}$ has leading
term $\bn^{\SB{k}}$.  Thus, $p_{\SB{k}}$ is a basis for $K[\bn]$, and
$q_{\SB{t}} = \sum_{\SB{k}} a_{\SB{k}} p_{\SB{k}}$. Thus, $Q_{\SB{t}}
= \sum _{\SB{k}} a_{\SB{k}} D^{\SB{k}} \bx^{\SB{t}+\SB{k}}$.
(Note that if $t_{i}<0$, then $a_{\SB{k}}$ is necessarily zero for
$k_{i}<-t_{i}$.)\qed 

Note that since $p_{\SB{k}}$ is a basis, the choice of $a_{\SB{k}}$ is
unique. Thus, we have the following result.
\begin{prop}[$\bX\bD$-Uniqueness] \label{coeffs}
XD-expansions and \bX\bD-expansions are unique. 
\end{prop}
\subsection{Closure of $DX$-Operators}

The main result of this section is the fact that the set $\DX $ of
operators with DX-expansions is 
closed under composition. 
The following corollary of proposition \ref{dxconv} is crucial to the
proof of theorem \ref{t4}. 
\begin{lemma}\label{lemon}
Using the above notation, if $Q\in \DX $, then
$q_{t}$ is identically zero for $t$ large.
\end{lemma}
{\em Proof:} Since $\lim_{k\rightarrow \infty} \ord(f_k)-k
= +\infty$, it follows that $k-\ord(f_k)$ takes a maximum value
$T$. Hence, 
$q_{t}$ is identically zero for $t>T$.\qed

It can similarly be shown that if $Q$  has a BX-expansion, then
$q_{t}$ is identically zero for $t$ large. 

\begin{thm}\label{t4}
The set $\DX $ of DX-operators 
 forms a (noncommutative, unitary)
$K$-subalgebra of the algebra ${\cal L}(K[x],K[x])$ of 
all linear operators equipped with the standard operations
\begin{eqnarray*}
(R \circ 
 P)p &=& R(Pp)\\
(cQ)p &=& c(Qp)\\
(R+P)p &=& Rp + Pp.
\end{eqnarray*}
\end{thm}
{\em Proof:} 
$\DX $ is trivially closed under multiplication by constants, and
under addition. We will now show that $\DX $ is closed under
composition. 

Let $Q$ and $R$ be DX-operators. Thus, as defined above $q_{t}$ and
$r_{s}$ are polynomials by theorem \ref{p3}. Moreover, by
lemma \ref{lemon}, $q_{t}$ is identically zero for 
$t>T$, and $r_{s}$ is identically zero for $s>S$. Define
$Q=R\circ 
 P$. It suffices to show that $q_{u}$ is a polynomial. 
\begin{eqnarray*}
Qx^{n} &=& \sum _{t\in {\DLZZ}} 
p_{t}(n)\, Rx^{n+t}\\
&=& \sum _{s,t\in {\DLZZ}} 
p_{t}(n)\, r_{s}(n+t)\, x^{n+t+s}\\
q_{u}(n) &=& \sum _{s+t=u} p_{t}(n)\, r_{s}(n+t)\\
 &=& \sum _{t=u-S}^{T} p_{t}(n)\, r_{u-t}(n+t).
\end{eqnarray*}
Thus, $q_{u}(n)\in K[n]$ since it is a finite sum of polynomials.$\Box$ 

Note that $q_{u}$ is identically zero for $u>S+T.$

Surprisingly, theorem \ref{t4} does not even generalize to two
variables. Consider the DX-operators $R=\sum _{k=0}^{\infty
}D_{1}^{k}X_{2}^{k} $, and $P=\sum _{k=0}^{\infty }
D_{2}^{k}X_{1}^{k}.$ Now, let $Q= R\circ P.$ Recall that $q_{00}(i,j)$
is the coefficient of $x_{1}^{i}x_{2}^{j}$ 
in $Qx_{1}^{i}x_{2}^{j}$. Thus,
$q_{00}(n,n)$ is given by the finite hypergeometric function
${}_{3}F_{0}(1,n+1,-n)=\sum _{k=0}^{n}(n)_{k}(n+k)_{k}$. 
All of the terms of the sum are nonnegative, and the term
corresponding to $k=0$ gives a lower bound of $(n!)^{2}$. This
guarantees
that $q_{00}$ is not polynomial. Thus, $Q$ is not a \bD \bX-operator 
by theorem~\ref{bolddxchar}.

\subsection{Coefficients of $DX$-Expansions}
Note that (\ref{star}) gives an efficient method of calculating the
XD-expansion of an operator $Q$. Given the action of $Q$  on a few
polynomials $Qx^{0},Qx^{1},\ldots,Qx^{n}$, the first few terms of its
XD-expansion can be automatically calculated \cite{us}. On the other
hand, the coefficients $f_{k}(D)$ of a DX-expansion (\ref{Dx}) depend
on the action of $Q$ on all the powers of $x$. No efficient means of
calculating DX-expansions is known in general. (Various special
techniques are used for each case treated in section \ref{s3}.) Is
there an analog of (\ref{star}) for DX-expansions? The only results in
this direction found so far are the following propositions.
\begin{prop}\label{new}
Let $Q\in \DX $. If
$$ Q=\sum _{n=0}^{\infty }c_{n}(D)\, X^{n}, $$
then 
$$ \frac{Q\exp(xt)}{\exp(xt)}=
\sum_{n=0}^{\infty}\sum_{k=0}^{n} {n\choose k}c_{n}^{(k)}(t)\,{x}^{n-k}.$$
\end{prop}
{\em Proof:} By (\ref{7.5}),
\begin{eqnarray*}
Q\exp(xt) &=& \sum _{n=0}^{\infty }c_{n}(D)\, X^{n} \exp(xt)\\
 &=& \sum _{n=0}^{\infty }\sum _{k=0}^{n }{n\choose k } X^{n-k}
c_{n}^{(k)}(D)\, \exp(xt)\\ 
 &=& \sum _{n=0}^{\infty }\sum _{k=0}^{n }{n\choose k } {x}^{n-k}
c_{n}^{(k)}(t)\, \exp(xt)\\ 
\hfill \frac{Q\exp(xt)}{\exp(xt)} &=& \sum _{n=0}^{\infty }\sum
_{k=0}^{\infty }{n\choose k } {x}^{n-k} c_{n}^{(k)}(t). \qquad\qed  
\end{eqnarray*}

By similar reasoning, we obtain the following proposition.
\begin{prop}\label{newer}
Let $Q\in \DX $.
$$ Q= \sum _{n=0}^{\infty}D^{n} a_{n}(X) ,$$
then
\begin{equation}
\frac{Q\exp(xt)}{\exp(xt)} =  \sum _{n=0}^{\infty }\sum _{k=0}^{n } {n \choose k}
a_{n}^{(k)}\, ({x})\, t^{n-k} 
\label{2usstar}. 
\end{equation}
\end{prop}
%
%
%
%
\section{Examples}\label{s3}

In this section we present example of DX- and XD-expansions.
\subsection{Umbral Operators}
A {\em delta operator} $P$ is a shift-invariant degree reducing operator.
Note that by the D-expansion formula $P=f(D)$ where $\ord (f)=1$.
Let $p_{n}(x)$ be the divided power sequence of $P$.
The operator $U_{P} \colon 
p_{n}(x)\rightarrow x^{n}$ is the {\em umbral
operator} associated to $P$, and the operator
$\sigma_{P}\colon 
p_{n}(x)\rightarrow (n+1)p_{n+1}(x)$  is the {\em umbral shift}
associated to $P$. For example, $U_{D}=I$ and $\sigma_{D}=X$.
\begin{prop}\label{p5}
The umbral operator $U_{P}$ is a DX-operator if and only if $Px=1$.
\end{prop}
{\em Proof:} {\bf (If)} Let $P=p(D)$, and let $r(t)=t-p(t)$. By
\cite[theorem 2.1]{GJ},  
$$U_{P} = \sum _{k=0}^{\infty } p'(D)\, r(D)^{k} X^{k}/k!.$$
By proposition \ref{dxconv}, the sum above converges since $\ord
(p'r^{k})=2k$ and $2k-k\rightarrow +\infty$.

{\bf (Only If)} Let $a=Px\in K$ and $Q=U_{P}$. The leading coefficient of
$p_{n}(x)$ is $1/a^{n}$. The result follows from theorem \ref{p3}
since $q_{0}(n)=a^{n}$ is not a
polynomial with respect to $n$ unless $a=1.$\qed 

On the other hand, all umbral operators have the following
XD-expansion.
\begin{prop}
Let $P$ be a delta operator. Then its umbral operator $U_{P}$ can be
expressed as
$$ U_{P} = \sum _{k=0}^{\infty } X^{k}(P-D)^{k}/k!. $$
\end{prop}
%
{\em Proof:} Let $\overline{p}_{k}(x)$ be the conjugate sequence of
polynomials for $P=p(D)$ defined by the generating function 
$$ \sum _{k=0}^{\infty }\overline{p}_{k}(x)\, t^{k}/k! = \exp(xp(t)). $$
Note by \cite[theorem~7, 
p.~708]{RKO}, $U_{P}x^{k}=\overline{p}_{k}(x)$. 
We now calculate the coefficients $c_{n}(x)$ of the XD-expansion by
applying (\ref{star}), 
\begin{eqnarray*}
\sum _{n=0}^{\infty } c_{n}(x)\, t^{n} &=& \frac{U_{P} \exp(xt)}{\exp(xt)} \\
&=& \exp(-xt)\, \sum _{k=0}^{\infty }(U_{P} x^{k})\, t^{k}/k! \\
&=& \exp(-xt)\, \sum _{k=0}^{\infty }\overline{p}_{k}(x)\, t^{k}/k! \\
&=& \exp(x(p(t)-t)). 
\end{eqnarray*}
Hence, $ U_{P} = \sum _{k=0}^{\infty } X^{k}(P-D)^{k}/k!$.\qed 

Contrast with the triple sum expansion given in \cite{Sun}. 
\subsection{Umbral Shifts}
\begin{prop}
All umbral shifts $\sigma_{P}$ are DX-operators. 
\end{prop}
{\em Proof:} Given $c\in K^{*}$, the umbral shift $\sigma_{P}$ is
identical to $c\sigma_{cP}$.
Thus, without loss of generality, we may suppose that
$Px=1$. Moreover, $\sigma_{P}= (U_{P})^{-1} X U_{P}$. Note that
$(U_{P})^{-1}= U_{R}$ where
$R=r(D),P=p(D),r(p(t))=p(r(t))=t$ (cf. \cite[theorem~7, 
p.~708]{RKO}).
By the Lagrange inversion formula, $Rx=1$. Thus, $\sigma_{P}$ is the
composition of three DX-operators: $U_{R}, X, U_{P}$. Hence, by
theorem \ref{t4}, $\sigma_{P}$ is a DX-operator.\qed 

{\em Alternate Proof:} The following XD-expansion of $\sigma_P$
is a restatement of Rodrigues' formula \cite[theorem~4, 
p.~695]{RKO}:
\begin{equation}\label{XDsig}
\sigma_P = X \frac{1}{P'}
\end{equation}
where $P'=P X - X P$ is the shift-invariant operator called
the  {\em Pincherle derivative} of $P$ \cite[Section~4, p.~694]{RKO}.
We then deduce
\begin{equation}\label{DXsig}
\sigma_P = \frac{1}{P'} X - \left( \frac{1}{P'} \right)'
= \frac{1}{P'} X + \frac{P''}{(P')^2} 
\end{equation}
as an explicit DX-expansion of $\sigma_P.$\qed 

Note that $X$ appears with exponent at most one in expansions
(\ref{XDsig}) and (\ref{DXsig}). 
\subsection{Endomorphisms}
\begin{prop}\label{p7}
The only endomorphisms of $K[x]$ with DX-expansions are the
translation operators.
\end{prop}
{\em Proof:} Let $Q$ be an endomorphism of $K[x]$, $Qp(x)= p(q(x)).$
If $\deg(q)>1,$ then 
$q_{t}$ is not identically zero for $t$ large, and by lemma
\ref{lemon}, $Q$ can not be a DX-operator. 

If $\deg(q)=0$, then without loss of generality (compose $Q$ with a
translation if necessary),  $Q$ is evaluation at zero. Thus,
$q_{0}(n)=1$ if $n=0$, and $0$ otherwise.  This function is not a
polynomial. Hence, by theorem \ref{p3}, $Q$ is not a DX-operator.

If $\deg(q)=1$, then $q(x)=ax+b$. Without loss of generality,
$q(x)=ax$. Thus, $q_{0}(n)=a^{n}$ which is not polynomial unless
$a=1.$\qed 

On the other hand, by (\ref{star}),  all endomorphisms $Q$ of $K[x]$
have the following  
XD-expansion
$$ Q = \sum _{k=0}^{\infty } (q(X)-X)^{k} D^{k}/k!.$$
where $Q(p(x))=p(q(x))$.
Notice the similarity to Taylor's formula.

In the multivariate case, let $q_{i}(x)=Qx_{i}.$ Then 
$$ Q = \sum _{\SB{k}} ({\bf q}(\bX)-\bX)^{\SB{k}}\,
\bD^{\SB{k}}/\bk !.$$ 
\section{Open Problems}\label{op}

\begin{description}
\item[a)] By theorem \ref{t4}, the product of two DX
operators is again a DX operator. Given the DX-expansion of $Q$ and
$R$, is there an explicit formula for the DX-expansion of their product?
Similarly, is there an explicit formula for the XD-expansion of any
two operators given their XD-expansions.

\item[b)]
If $K$ is a topological field, then we have a weaker notion of
convergence than the discrete topology, and thus more DX-operators.
How can they be classified?

\item[c)]
Is the product of two BX-operators also a BX-operator? This problem is
especially difficult since we have no criteria analogous to theorem
\ref{p3} to tell whether a product is a BX-operator. In fact, it
remains to be seen whether BX-expansions are unique (when they exist).

\item[d)]
Find a formula more explicit than propositions \ref{new} and
\ref{newer} by which to calculate DX-expansions. Is there a BX-analog of 
propositions \ref{new} and
\ref{newer}?

\item[e)]
Let $B$ be a degree lowering operator as above, and let $Y$ be a
degree raising operator. For what $B$ and $Y$ can all linear operators
be expressed by a (unique) YB-expansion: $\sum_{k=0}^\infty B^k
a_k(Y).$    

\item[f)] Characterize operators $Q$ with identical XD and DX-expansions.
That is, 
$$ Q=\sum _{i,j}c_{i,j}X^{i}D^{j}= \sum _{i,j}c_{i,j}D^{j}X^{i}.$$
Conjecture:  $Q$ must be of the form $p(X)+f(D)$ for some
polynomial $p$ and formal power series $f$.

\item[g)] Extensions of the Umbral Calculus from polynomials to inverse
Laurent series (negative powers of $x$) and Artinian series
(fractional powers of $x$) are known \cite{Roman,Loeb}. Can these be
used to derive $XD$ expansions of operators that act not on
polynomials, but on Laurent or Artinian series?
\end{description}

\end{document}